\documentclass[letterpaper, 10 pt, conference]{ieeeconf}

\IEEEoverridecommandlockouts   
\usepackage{amsmath}
\usepackage{subfig}
\usepackage{amssymb}
\usepackage{graphicx}
\usepackage{color}
\usepackage{url}
\usepackage{bbm}
\usepackage{algorithm2e}

\newcommand{\bd}{\mathbf}

\newcommand{\E}{\mathbb{E}}

\title{\LARGE \bf Mitigation of Coincident Peak Charges via Approximate Dynamic Programming
}

\author{Chase P. Dowling and Baosen Zhang
\thanks{Department of Electrical and Computer Engineering, University of Washington Seattle, WA 98195, USA	 \{cdowling, zhangbao\}@uw.edu}%
\thanks{This work is partially supported by grants from Centrica, plc.}%
}

\begin{document}
\maketitle
\thispagestyle{empty}
\pagestyle{empty}

\begin{abstract}
A significant portion of a consumer's annual electrical costs can be made up of coincident peak charges: a transmission surcharge for power consumed when the entire system is at peak demand. This charge occurs only a few times annually, but with per-MW prices orders of magnitudes higher than non-peak times. While predicting the moment of peak demand charges over the course of the entire billing period is possible, optimal cost mitigation strategies based on these predictions have not been explored. In this paper we cast coincident peak cost mitigation as an optimization problem and analyze conditions for optimal and near-optimal policies for mitigation. For small consumers we use approximate dynamic programming to first show the existence of a near-optimal policy and second train a neural policy for curtailing coincident peak charges when subject to ramping constraints.
\end{abstract}

\section{Introduction}
\label{sec:intro}

A coincident peak (CP) is a consumer's electrical demand at the time of the total system peak demand. Since much of the power system infrastructure is only used only during peak times~\cite{uddin2018review}, some system operators and utilities use CP pricing mechanisms to incentivize customers to reduce their consumption during peak times, therefore hoping to achieve an overall reduction of the system peak~\cite{ercot4cp, fortcollinscp}. Existing CP charges are applied through a rate structure, with the rates at peak times hundreds of times larger than at regular times. As a result, CP charges often account for a significant portion---often greater than 20\%---of annual electrical costs for participating customers~\cite{liu2013data}, providing them with a strong incentive to reduce their consumption at these peak times~\cite{consumerresp, zarnikau2013response}. 

In this paper, we adopt the view point of a small customer facing CP charges and study how the customer can operationally mitigate this cost. The primary challenge is that the timing of the CP charges are only known after all of the system demands have been realized. For example, if CP is charged on a monthly basis~\cite{fortcollinscp}, the hour that the peak load occurred in only determined after the entire month has passed. 

To mitigate this uncertainty in peak timing, operators typically provide warning signals to consumers to indicate peak is forthcoming. In~\cite{liu2013data}, the authors utilize these signals to develop a scheduling model for a data center's workload in the Fort Collins PUD \cite{fortcollinscp}. However, forecasting when a peak will occur is a difficult prediction problem~\cite{ercotloadforecast, dowling2018coincident}, since it only occurs (by definition) at a single point in time. Since the rate associated with the CP is orders of magnitude higher than normal time-of-use rates, false negative predictions are extremely costly. Therefore operators tend to send out many successive CP warning signals, degrading the efficiency of customer responses and leading to user fatigue in the long run~\cite{zarnikau2007industrial,zarnikau2013response}.

In this work we treat the problem of mitigating CP costs as an optimization problem that is continually solved over the entire horizon of the billing period. Instead of explicitly predicting when the peak will occur, we adopt a probabilistic framework to gracefully incorporate observations made by the customer to maximize their expected revenue. That is, at each time-step we calculate the probability of the peak occurring at some point in the future having observed previous values of system demand. 

Related works on mitigating CP pricing focus on large consumers with considerable demand flexibility, namely, data centers~\cite{wierman2014opportunities}. Limited works have addressed CP prices for data center consumers directly such as \cite{liu2013data} which incorporates existing grid operator signals. Others related to data center peak power consumption address the problem generally based on time-of-use costs given on-site storage or generation capabilities \cite{shi2016leveraging, wang2013data} without tackling the idiosyncrasies of CP pricing mechanisms. 

Dynamic programming is a natural approach to maximize the expected revenue of a small customer in the face of CP timing uncertainty. However, since the action space of a customer is continuous and coupled in time, solving the dynamical programming problem becomes intractable. Therefore we approximate the value function and train a deterministic policy parametrized as a neural network. Based on the structure of the CP charge, we design the input of the neural network to explicitly include the maximum of the observed demand and the number of time periods. Using these inputs, we show that this neural network based policy is comparable to a brute force grid search and outperforms a standard benchmark algorithm. This approach advances the state-of-the-art by providing a way to actively reduce the CP cost that does not rely on system warning signals or assumes an adversarial environment.  


The rest of the paper is organized as follows: Section~\ref{sec:problem} defines the optimization problems to be solved, Section~\ref{sec:model} provides the solution framework, Section~\ref{sec:results} presents a numerical case study. We conclude with a discussion on future work for both large and small consumers and make some final remarks in Section~\ref{sec:conclusion}.

\section{Model and Problem Formulation}
\label{sec:problem}
We consider a small customer that tries to maximize its revenue subject to CP charges over $T$ time periods. Let $x_t$ be the energy consumption of the customer at time $t\in\{1,\dots,T\}$, and it is limited to be between $\underline{x}$ and $\overline{x}$. We assume the revenue of the customer is represented by a concave increasing function $g(\cdot)$~\cite{kirschen2018fundamentals,ratliff2014incentive} of power consumption. Let $\pi_{cp}$ be the CP charge rate and the customer pays an amount of $\pi_{cp} x_{t^*}$ where $t^*$ is the time period the system peak occurs. 

\begin{figure}
    \centering
    \subfloat[]{\label{fig:forecast}\includegraphics[width=0.9\linewidth]{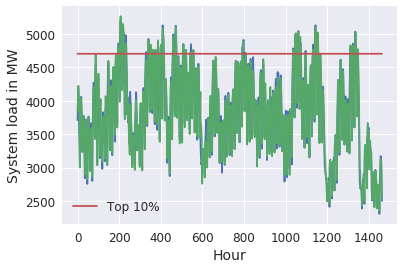}} \\
    \subfloat[]{\label{fig:errorhist}\includegraphics[width=0.9\linewidth]{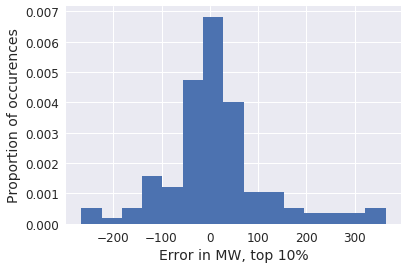}}
    \caption{In PJM's Duke Energy Ohio/Kentucky region: (a) actual and forecasted system load and (b) the distribution of forecast error from June 1 to September 30, 2018.}
    \label{fig:noisedist}
\end{figure}

We model the system load in each time period as random variables $S_1,\ldots ,S_T$ where the mean of $S_t$ is the forecasted load value. Even though system loads are strongly correlated in time, once the forecast value is given, the forecast errors are typically independent across time periods~\cite{feinberg2005load,weron2007modeling}. For example, Fig.~\ref{fig:noisedist} illustrates the distribution of system load forecast error for the top 10\% of system load values during summer months in 2018 in a single PJM subregion.

Additionally, it is clear a large customer has a measurable impact on the value of $S_t$. In the case of ERCOT, for example, the 6 largest customers accounted for over 80\% of each summer monthly CP in 2017 and 2018. On the other hand, of the 130 total customers participating in ERCOT's CP pricing program, 20\% in 2017, and 40\% in 2018 consumed less than $\frac{1}{10}$ the difference between the annual system peak and the second closest system demand. For these exceedingly small customers, the variance of the forecast error well exceeds their CP demands of less than 5-10 MW.

Therefore, in the case of a small customer we assume that $S_1,\dots,S_T$ are independent and their mean is given by the forecast values. We define $t^*$ as the time index corresponding to the maximum load:
\begin{equation}\label{eqn:ts}
t^*=\arg\max_t \{ S_1,\dots,S_t\}. 
\end{equation}
Note since $t^*$ is a function of random variables, it is also a (discrete) random variable. 

With these definitions, the expected net revenue (or reward) of a customer is
\begin{equation}\label{eqn:reward}
    \E[R] :=\E\left[\sum_{t=1}^{T} g(x_{t}) - \pi_{cp}x_{t^{*}}\right]
\end{equation}
where $t^*$ is defined in \eqref{eqn:ts} and the expectation is over the random variables $S_1,\dots,S_T$. 

The goal of the customer is then to maximize $\E[R]$ subject to their operational constraints. We assume a fairly simple customer model where the demand from each time-period is coupled through a ramping constraint, and the customer's optimization problem is
\begin{equation}\label{eqn:optprob1}
\begin{aligned}
& \underset{x_{t}}{\text{maximize}}
& & \E[R] \\ 
& \text{subject to}
& & \underline{x} \leq x_{t} \leq \bar{x}\\
& & & x_{t-1} - \delta \leq x_{t} \leq x_{t-1} + \delta
\end{aligned}
\end{equation}
where $\delta$ limits the possible rate-of-change between two time periods. Other types of constraints can be included using the techniques described in this paper. 

\subsection{Sequential Decision Problem}
In practice, the optimization problem in \eqref{eqn:optprob1} needs to be solved in a sequential manner. There are two sources of temporal coupling in \eqref{eqn:optprob1} that makes this sequential optimization problem nontrivial and interesting. The first is the ramp constraint between two successive time steps. The second is how the timing of the peak changes after loads are observed. 

At time $t$, the customer have observed the realization of $S_1,\dots,S_{t-1}$, which we denote as $s_1,\dots,s_{t-1}$. Based on these observations, the value of peak time $t^*$ changes. In other words, if the observed loads are large, then the peak is likely to have already occurred and the customer can act aggressively; conversely, if the maximum value of the observed loads are small, then the customer should act more conservatively to protect against incurring a large CP charge in the future. Therefore, even if the ramp limits are not present, the structure of the CP charge induces a time dependency on the observations made at each stage. A variant of \eqref{eqn:optprob1} without ramp constraints is studied in \cite{liu2013data} where data centers are assumed to have a large amount of flexibility and can ignore the coupling of its own actions between two time periods. In this paper, we focus on commercial customers that may not have this level of flexibility and are limited in their rate-of-change.

\subsection{Benchmark Algorithm}
There are two typical strategies to solve \eqref{eqn:optprob1} in practice. The first is to simply assume that all time periods (e.g., all hours between 3 PM and 7 PM on a hot summer day) experience the peak demand and conservatively reduce the load to mitigate the CP charge~\cite{zarnikau2013response}. The the CP charge is evenly distributed over all of these time intervals. The second is to follow the warning signals of operators and treat those as true peak times~\cite{fortcollinscp,liu2013data}. It turns out that these two strategies amount to the same thing, since operators tend to be conservative and issue CP warnings for all of the time periods that have a reasonable chance of experiencing the moment of peak demand~\cite{fortcollinscp}. This is to say that conservative CP warnings amount to treating any hot summer afternoon, for example, as equally likely to the system peak without taking into joint consideration the known system capacity and previously observed system loads during the billing window. Therefore we adopt the following strategy as the a baseline algorithm which we call the naive strategy \cite{dowling2018coincident}, where the customer solves
\begin{equation} 
 \max_{x_t} \; g(x_t)- \frac{1}{T} \pi_{cp} x_t,
 \end{equation} 
 where the scaling factor $1/T$ represents the fact that the cost of CP is amortized evenly to all of the time periods under consideration. The optimal solution is then the the demand that satisfies the first order optimality condition $Tg'(x^*)-\pi_{cp}=0$. 
 
 Note even though this solution is simple to compute, it does not take into account the successive realization in the system load and is generally suboptimal. In later comparisons, we will call it the naive policy. In the next section, we develop a policy based on approximate dynamic programming to solve \eqref{eqn:optprob1}.

\section{Approximate Dynamic Programming}
\label{sec:model}
\subsection{Dynamic Programming Formulation}

Let us first directly apply a dynamic programming approach to optimize \eqref{eqn:reward}. Suppose the customer is solving for the optimal $x_{T}$, having already chosen $x_{1},\ldots, x_{T-1}$ at the final step $t = T-1$. The customer must maximize the expected reward conditioned on observed system load realizations, $s_1, \ldots, s_{T-1}$, specifically, $\mathbb{E}[R|s_{1},\ldots, s_{T-1}]$.

At $t = T-1$,  let $s_{m} = \max\{s_1,\ldots, s_{T-1}\}$, the maximum observed so far; since this is the final round, the expected reward depends \emph{only} on whether $s_{T}$ will be larger than $s_{m}$. Let $p_{T} = 1 - P(s_{m} < S_{T})$, the probability that the final system load realization will be the CP. Then the objective,

\begin{subequations}
\begin{align}
     \mathbb{E}[R|s_m] &= \sum_{t = 1, \ldots T - 1} g(x_{t}) + g(x_{T}) - \pi_{cp}\mathbb{E}\left[ x_{t^{*}} | s_{m} \right] \\
    &= \sum_{t = 1, \ldots T - 1} g(x_{t}) + g(x_{T}) - \\
     &\qquad \qquad \qquad \pi_{cp} [ (1 - p_{T})x_{t^{*}} + p_{T}x_{T}]. \label{eqn:dynfinal}
\end{align}
\end{subequations}

Thus, for the solution $x'_{T}$ to $g(x'_{T}) - \pi_{cp}p_{T} = 0$, the optimal $x^{*}_{T}$ is the point in the interval $[x_{T-1} - \delta, x_{T-1} + \delta]$ which minimizes $|x'_{T} - x^{*}_{T}|$. At $t = T-2$, in order to solve for the optimal $x^{*}_{T-1}$ there are two potential rounds that the CP may yet occur on and the customer must consider the probability that either $S_{T}$ or $S_{T-1}$ is the CP. Indeed, 

\begin{subequations}
\begin{align}
    \mathbb{E}[R|s_{1},\ldots, s_{T-2}] &= \sum_{t=1,\ldots,T-2} g(x_{t}) + g(x_{T-1}) + \\
    &\mathbb{E}[g(x_{T}) - \pi_{cp}[(1 - p_{T})x_{t^{*}} + p_{T}x_{T}]],
\end{align}
\end{subequations}

\noindent noting that $x_{T}$ remains inside the expectation since it depends on the realization of $S_{T-1}$. Iterating backwards yields a dependency on future realizations of $S_{t}$, where only the current consumption $x_{t}$, maximum system load observed thus far $s_m$, and number of rounds remaining $T-t$ influence future choices of $x_{t+1},\ldots, x_{T}$. 

A straightforward means of addressing this would be a brute force grid search. Consumption values in $[\underline{x}, \bar{x}]$ and a range of likely system loads $S_{t}$ can be discretized, with every potential outcome being computed forward from each possible initial value $x_{1}$. At each time a consumer would choose a feasible $x_{t+1}$ subject to ramping constraints that maximizes the expected reward over the \emph{entire} horizon $T$ for all possible outcomes given $s_{1},\ldots, s_{t}$ using the output of the grid search as a look-up table; however, a complete grid search exhibits exponential complexity in $T$. This dimensionality problem is common in applications of dynamic programming \cite{si2004handbook}.

Therefore we propose an approximate dynamic programming approach by sampling from all possible outcomes in order to estimate the best choices of $x_{t+1}$. These samples are used to train a policy $f$, which takes as input at time $t$ the current consumption $x_{t}$, the largest system load observed so far in the billing period $s_{m}$, and the number of rounds left, $T - t$. The policy then outputs an estimated optimal $\hat{x}_{t+1}$. We note that in the absence of these time coupling constraints, an optimal solution exists since each time-step is completely independent.

\subsection{Neural Network Policy}
Neural networks have gained popularity as a tractable way to parameterize policies. For example, they have been used to solve approximate dynamic programming problems in \cite{si2004handbook, bertsekas1996neuro, wang2009adaptive}. We also adopt a neural network based policy to solve \eqref{eqn:optprob1}. In the context of dynamic programming, a policy is a function that maps previous values to an action. In our case, this policy should map the current choice of $x_{t}$ and observations of $s_{t}$ to an output $x_{t+1}$ that a customer should select as their demand based on, in this case, criterion that maximizes their expected reward over the remaining time horizon. A policy in the context of approximate dynamic programming attempts to output a value $\hat{x}_{t+1}$ that is close to optimal.

\begin{figure}
\vspace{1em}
    \centering
    \includegraphics[width=0.8\linewidth]{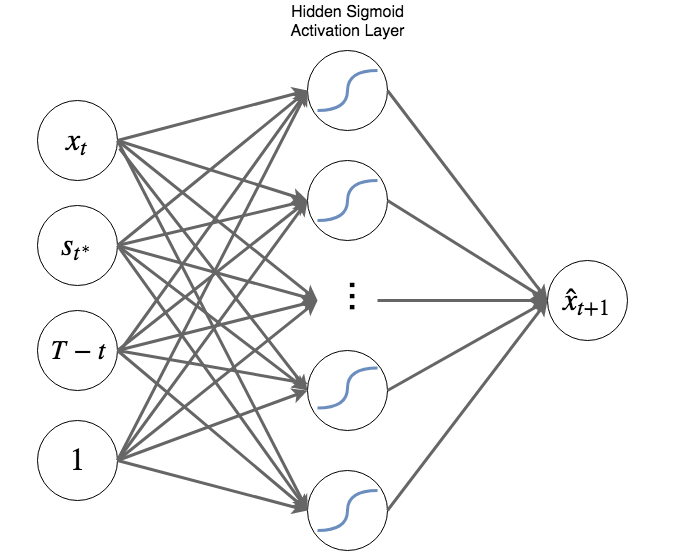}
    \caption{Architecture of single-layer neural network policy, with inputs $x_{t}$, $s_{m}$, $T - t$, and a linear bias term.}
    \label{fig:nnf}
\end{figure}

Therefore, in order to train a policy $f$ we require an approximation of the true optimal output $x^{*}_{t+1}$ of $f$ that maximizes expected reward $R$ given previous observations of $s$. Alg.~\ref{alg:samp} details the process by which these samples are generated. At time $t$, for each feasible value of $x_{t+1}$ subject to the ramping constraints, potential outcomes are forward simulated until time $T$ a total of $C$ times. The value of $x_{t+1}$ with the best average \emph{remaining} reward---modulo $\sum_{i = 1}^{t} g(x_{t})$---is selected as the training output.  If the range of customer consumption values are discretized into $n$ values, sampling across all possible starting times $t \in T$ yields an improved complexity of $\mathcal{O}(TCn)$.

Many samples are compiled and used to train a neural network to to approximate the function $f$. Fig.~\ref{fig:nnf} illustrates the basic network architecture described and used in Sec.~\ref{sec:results}. The necessary inputs of the policy at time $t$ are 1) the current state $x_{t}$, as this determines the range of feasible values due to the ramping constraint, 2) the maximum value $s_m = \max\{s_1, \ldots, s_t\}$ observed thus far, as values of $s_i < s_m$ have no bearing on the timing of the CP, and 3) the number of rounds left, $T - t$---given the probability density function of $S_{t}$---determines the probability any future value will be greater than $s_m$ and thus the new potential CP. When performing grid search, these three values completely determine the expected reward when choosing a value $x_{t+1}$.

\begin{algorithm}[t!]
\SetAlgoLined
\KwData{$x_{t}, s_{m}, T$}
\KwResult{Estimated policy output $\hat{x}_{t+1}$}
 $C = $ number of Monte Carlo simulations\;
 sim\_rewards $= []$\;
 discretize $[x_{t} - \delta, x_{t} + \delta]$\;
 \For{each feasible $x_{t+1} \in [x_{t} - \delta, x_{t} + \delta]$}{
  sim\_rewards$[x_{t+1}]$ $\leftarrow []$\;
   \For{$j = 1,\ldots, C$}{
     sample $s_{t+1}$ according to system load distribution\;
     $\bd{x} \leftarrow [x_{t}, x_{t+1}]$\;
     $\bd{s} \leftarrow [s_{m}, s_{t+1}]$\;
     \For{$k = t+2,\ldots, T$}{
     sample $s_k$ according to system load distribution\;
     $\bd{s} \leftarrow s_{k}$\;
     randomly sample feasible $x_{k}$ from interval $[x_{k-1} - \delta, x_{k-1} + \delta]$\;
     $\bd{x} \leftarrow x_{k}$
     }
     sim\_rewards$[x_{t+1}] \leftarrow R(\bd{x}, \bd{s})$\;
    }
   }
   $\hat{x}_{t+1} = \arg\max_{x_{t+1}} \frac{1}{C} \sum$ sim\_rewards (choose $x_{t+1}$ with best average forward simulated reward)
\caption{Monte Carlo path sampling}\label{alg:samp}
\end{algorithm}

\section{Case Studies}
\label{sec:results}

\begin{figure}
    \centering
    \includegraphics[width=0.9\linewidth]{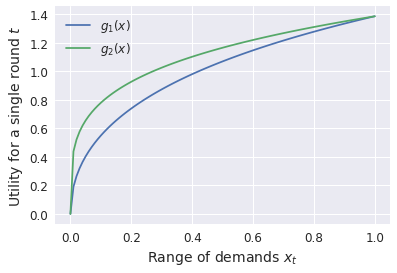}
    \caption{Case study revenue functions $g_{i}(x)$}
    \label{fig:utilfuncs}
\end{figure}

To test the efficacy of our approximate dynamic programming solution compared to the naive strategy, we set up a numerical study\footnote{All code to reproduce our results and results plots can be found on our GitHub project repo at \url{https://github.com/cpatdowling/peakload/blob/master/notebooks/cdc_2019_submission.ipynb}}. We consider two different consumer revenue functions (illustrated in Fig.~\ref{fig:utilfuncs}),

\begin{subequations}
\begin{align}
    g_{1}(x) &= 2\log(1 + x^{2})\,\,\, \textnormal{and}\\
    g_{2}(x) &= 1.386 \sqrt[\leftroot{-2}\uproot{2}4]{x}
\end{align}
\end{subequations}


For both revenue functions we suppose the customer's ramp constraint $\delta = 0.3$, and that the customer's CP charge rate $\pi_{cp} = 0.6 T g(\bar{x})$, or 60\% of their maximum possible gross revenue over $T$ rounds. Typically CP charges form greater than 20\% of their annual electrical costs, but we choose to much higher percentage to illustrate a more drastic scenario.

For $T = 2,... 10$ rounds\footnote{Results for larger values of $T$ can be found for $g_{1}(x)$ in our Github project repo}, we use the sampling strategy defined in Alg.~\ref{alg:samp} to generate 1000 input/output samples per time $t \in [1,\ldots, T]$, such that we train with an even number of $\hat{x}_{t+1}$ for all $t$. For each feasible $x_{t+1}$ being evaluated, the number of simulations $C = 100$. 

\begin{figure}
    \centering
    \includegraphics[width=0.9\linewidth]{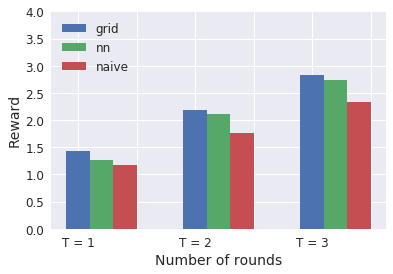}
    \caption{Comparison of best-possible performance via grid search to a NN policy and the naive strategy. Reward is strictly increasing with each additional number of rounds roughly as $Tg_{1}(x)$}
    \label{fig:perfcomp}
\end{figure}

%

As a first pass we design our neural network to have a single hidden layer of size $4$ with a sigmoidal activation function. Further, we included a linear bias term, indicated as the fourth input in Fig.~\ref{fig:nnf}. The neural network is trained using mean-squared error; additional hyperparameters like learning rate and batch-size can be found in our linked repository.

First we test the validity of the assumption that our ADP sampling procedure yields near-optimal choices $\hat{x}_{t+1}$ against an exhaustive grid-search. For $T = 2, 3$ and $4$, and revenue function $g_{1}(x)$ we compute an exhaustive grid for our example function and system load distribution. Fig.~\ref{fig:perfcomp} illustrates the relative performances of each strategy; while we found that the discretization resolution of the grid search has a noticeable effect on the resulting reward given the choice of our case study revenue function, the NN policy performs nearly as well and we make use of the sampling technique on a larger number of rounds for which grid search is intractable.

\begin{figure}
    \centering
    \includegraphics[width=0.9\linewidth]{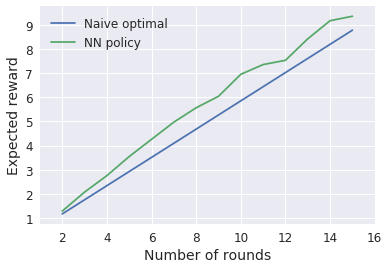}
    \caption{Policy performance for revenue $g_{1}(x)$, across time horizons $T$ with $\pi_{cp}$ set to be 60\% of maximum, unpenalized revenue for each $T$.}
    \label{fig:policyperf}
\end{figure}

\begin{figure}
    \centering
    \includegraphics[width=0.9\linewidth]{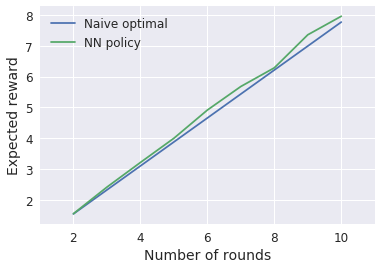}
    \caption{Policy performance for revenue $g_{2}(x)$, across time horizons $T$ with $\pi_{cp}$ set to be 60\% of maximum, unpenalized revenue for each $T$.}
    \label{fig:policyperfpoly}
\end{figure}

Figures \ref{fig:policyperf} and \ref{fig:policyperfpoly} illustrate the performance of the respective NN policies against the naive strategies for $g_{1}(x)$ ad $g_{2}(x)$. The NN policies consistently outperforms the naive optimal solution while maintaining the added benefit of being an solution approximated from sampled paths. In the case of mitigating CP costs on an hourly basis, it is unlikely that a potential CP would occur at anytime in excess of 8 to 10 consecutive hours, viz occuring outside known, afternoon peak hours \cite{liu2013data}. This benefits a customer by allowing them to focus on sweeping training parameters to tune the policy for a narrow range of values of $T$.

\begin{figure}
    \centering
    \includegraphics[width=0.9\linewidth]{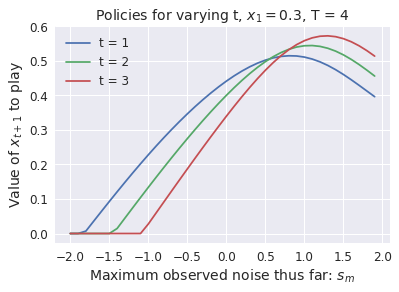}
    \caption{Example of policy for revenue $g_{1}(x)$ for multiple rounds $t$ over time horizon $T = 4$ and fixed $x_t = 0.3$. With later rounds of $t$, the policy becomes less conservative and shifts to the right as the decreasing number of rounds decreases the probability of a new maximum system load being observed.}
    \label{fig:examppolicy}
\end{figure}

\begin{figure}
    \centering
    \includegraphics[width=0.9\linewidth]{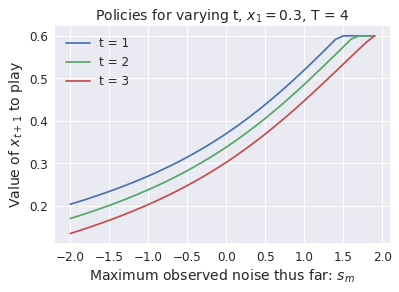}
    \caption{Example of policy for revenue $g_{2}(x)$ for multiple rounds $t$ over time horizon $T = 4$ and fixed $x_t = 0.3$. With later rounds of $t$, the policy interestingly becomes \emph{more} conservative, likely due to the sharper decrease in $g'_{2}(x)$ in increasing $x$.}
    \label{fig:examppolicypoly}
\end{figure}

Figures \ref{fig:examppolicy} and \ref{fig:examppolicypoly} illustrate examples of the outputs of $f$ corresponding to $g_{1}(x)$ and $g_{2}(x)$ both trained for $T = 4$. For each time $t = 1, 2, 3$, the output $x_{t+1}$ is given as a function of $s_m$ for a fixed initial value, $x_t = 0.3$. Note that with each consecutive round, the likelihood of $S_{t}$ remaining that may be a CP changes both as a function $t$ \emph{and} $s_m$. In general, the probability that the next realization of $S_{t}$ will be the maximum over all $T$, $p_{t} = \frac{1}{T - t}(1 - P(S_{t} < s_m)^{(T - t)}$. Interestingly these policies learned for $g_{1}(x)$ and $g_{2}(x)$ appear very different.

In the case of policy $f$ learned for $g_{1}(x)$, outputs of the policy for each $t$ become \emph{less} conservative with decreasing number of rounds. The flattening of the policy at small values of $s_m$ is due to the ramping constraint. Conversely, decreasing the radius of curvature of the revenue function---the sharpness in the initial revenue increase transitioning into diminishing returns---as in $g_{2}(x)$, it appears that a \emph{more} conservative strategy arises, likely due to the decreased cost of false negatives. Large changes in $x$ result in little change in $g_{2}(x)$ but correspondingly a relatively larger marginal decrease of $\pi_{cp}$ in the CP charge. That is to say it costs the customer little to curtail to values already near the naive optimal strategy for $g_{2}(x)$, (e.g. for $T=10$, $x = 0.311$), yet the NN policy still improves on the naive strategy.

\section{Conclusion}
\label{sec:conclusion}
In sum we considered how a small customer participating in a CP pricing program can near-optimally trade off lost revenue for CP cost savings. We formulated an approximate dynamic programming problem that incorporates successive observations of system loads likely to be a CP as inputs allowing a customer subject to ramping constraints to make informed curtailment decisions over the course of the billing period time horizon. This work improves on existing algorithms such as \cite{liu2013data} or \cite{dowling2018coincident} by escaping an ad-hoc threshold curtailment regime where if some measure exceeds a threshold parameter than the customer should curtail. Further, this optimization framework provides footing for further theoretic analysis, detailed below.

\subsection{Future Work}

The goal of future work is to explore the implications of the optimization formulation in \eqref{eqn:optprob1} for large players as observations of $S_t$ become a function of $x_{t}$. 
Multiple large customers contributing to $S_t$ resemble a Cournot competition, and a desirable outcome might be the existence of convergent strategies for large customers.

Additionally, for both large and small players, given a customer's time-coupled revenue function, how are they incentivized to participate in a CP program? Intuitively a customer with more demand flexibility---such as in the case of the data center in \cite{liu2013data}---has potentially more to gain from participating in a CP pricing program than a comparably sized customer with little demand flexibility. If this is the case, how do factors such as CP billing horizon affect incentives, e.g. annually vs. monthly, impact this incentives like discounted time of use rates in exchange for participating?

\bibliography{main}
\bibliographystyle{IEEEbib}

\end{document}